\documentstyle[12pt]{article}
\begin{document}
\begin{titlepage}
\rightline{Helsinki Institute of Physics preprint,  HIP-1997-65/TH}
\rightline{November 1997}
\rightline{hep-th/9712xxx}
\vskip 1.5cm

\begin{center}
{\huge \bf Dynamical systems related to }\\[.5cm]
{\huge \bf the Cremmer-Gervais $R$-matrix}
\end{center}
\vskip .6cm

\centerline{\Large \bf M.Chaichian${}^{1,2}$, 
P.P.Kulish${}^{2}$\renewcommand{\thefootnote}{*}\footnote{On leave
of absence from the St.Petersburg 
Department of the Steklov
Mathematical Institute, 
Fontanka 27, St.Petersburg, 191011,
Russia; ( kulish@pdmi.ras.ru )} }
\vskip .3cm
\begin{center}
${}^1$ High Energy Physics Division, Department of Physics\\ 
${}^2$ {Helsinki Institute of Physics}\\
{FIN-00014, University of Helsinki, Finland}
\end{center}
\vskip .3cm
\centerline{\bf and}
\vskip .3cm
\centerline{\Large \bf E.V.Damaskinsky}
\vskip .3cm
\begin{center}
{University of Defense Constructing 
Engineering, St.Petersburg, Russia}
\end{center}
\vskip 1.0cm

\begin{abstract}
The generalized Cremmer-Gervais $R$-matrix being a twist of 
the standard  $R$-matrix of $SL_q(3)$, depends on two extra parameters. 
Properties of this $R$-matrix are discussed and two dynamical systems, 
the quantum group covariant $q$-oscillator and an integrable spin 
chain with a non-hermitian Hamiltonian, are constructed. 
\end{abstract}

\end{titlepage}

\section{Introduction}
\renewcommand{\theequation}{\thesection.\arabic{equation}}
\setcounter{equation}{0}

A variety of integrable models of the quantum field theory,
classical statistical mechanics and quantum mechanics (systems 
with finite degrees of freedom) are studied using the quantum 
inverse scattering method (QISM) (see e.g. \cite{LDF}-\cite{Iz}).  
The main ingredient or a cornerstone of their
solution in the framework of the QISM is the $R$-matrix, a
solution to the Yang-Baxter equation (YBE). Although there is 
plenty of known solutions, in particular, corresponding to any
simple Lie (super)algebra \cite{KulSk}, \cite{YBE}, \cite{DrQG},
any new $R$-matrix, even being connected one way or another to
the known ones gives rise to a new set of integrable models. A
simple method to get ''new'' $R$-matrices from the standard ones
(constant or spectral-parameter- dependent) is the twist procedure
\cite{DR-Tw}, \cite{ReshTw}. The multiparametric $R$-matrix  
corresponding to the algebra $gl(n)$ with $\frac12 n(n-1)$ extra parameters 
on the diagonal, was obtained from the standard $sl_q(n)$ $R$-matrix 
( $\omega = q - q^{-1}$)
\begin{equation}
R(q)= \sum\limits_{i=1}^n  q\,e_{ii}\otimes e_{ii} + 
\sum\limits_{i=j}  e_{ii}\otimes e_{jj} + 
\sum\limits_{i < j}  \omega \,e_{ij}\otimes e_{ji} \;, 
\label{a1.1}
\end{equation} 
where $(e_{ij})_{ab} = \delta_{ia} \delta_{jb}$ are basis matrices, 
by the twist $F$ constructed from the generators 
of the Cartan subalgebra \cite{ReshTw}
\[
R(q,\{p_{ij}\})=F_{21}R(q)F_{12}^{\quad -1},\qquad
F\in {\cal U}({h})\otimes {\cal U}({h}).
\]
In this paper we shall look for integrable systems 
related to a particular $R$-matrix, so we shall
write formulas in some fixed irreducible representations and not
in the universal algebraic form of the quantum group theory
\cite{DrQG}, \cite{FRT}.

Recently it was shown that the $R$-matrix $R_{CG}(q)$ derived for the
quantum Toda theory \cite{Cr-Ger} is a twist of the 
standard $R$-matrix of the $sl_q(n)$ \cite {Hodge}, \cite {Jacob}. 
We shall restrict 
ourselves to the case of the $sl(3)$ algebra. 

An interesting feature of this generalized Cremmer-Gervais
$R$-matrix $R(q,p,\nu)$ (the twist adds two more parameters to the
original solution \cite{Cr-Ger} of the YBE) consists of two new
non-zero entries with respect to the standard $R(q)$
\begin{equation}
\begin{array}{ll}
R(q,p,v)= & R(q)+(p-1)\left( e_{11}\otimes e_{22}+e_{22}\otimes
e_{33}\right) \\
& +(p^{-1}-1)\left( e_{22}\otimes e_{11}+e_{33}\otimes
e_{22}\right) + ( p^2 /q - 1)e_{11}\otimes e_{33}\\
& + ( q/p^2 -1)e_{33}\otimes
e_{11}+q \nu \,(e_{32}\otimes e_{12}- p^2 /q^2 e_{12}\otimes e_{32} ) \,.
\end{array}
\label{f1}
\end{equation}
Due to these extra entries proportional to $\nu$ 
the quantum space ${\bf C}_q^3$
commutation relations for ($\ref{f1}$), which were 
obtained already in \cite {Cr-Ger},  
are more complicated than
for the $sl(3)$ case: $ x_ix_j=q_{ij}x_jx_i,$ $1\leq i<j\leq 3.$
They coincide with a deformed Heisenberg algebra 
( $q$-oscillator ) ${\cal A}_q$ 
\begin{equation} 
x_2x_1=pqx_1x_2,\quad
x_2x_3=(pq)^{-1}x_3x_2,\quad x_1x_3-p^{-2}x_3x_1=\nu\,x_2^{\;2} 
\label{f2} 
\end{equation} 
after the identification $x_1 = A\,, x_2 = K \,, x_3 = A^{\dagger}$ 
with some set of ${\cal A}_q$ generators (Sec.3). 
Looking at a quantum group covariance property of this system:
$x_i\rightarrow \sum_jT_{ij}\otimes x_j$ and studying the
corresponding differential calculus \cite {WZ} the $R$-matrix (\ref{f1})
was found also in \cite{Irac}.

There is no quadratic element (Hamiltonian) in ${\cal A}_q$ 
invariant with respect to the quantum group ${\cal A}(R)$ 
coaction:  $\varphi (x_i)=\sum_jT_{ij} \otimes x_j$. Hence, 
after coaction one gets more complicated Hamiltonian in 
a tensor product of representation spaces of ${\cal A}(R)$ and 
${\cal A}_q$. We analyze in this paper 
the physical consequences of this transformation for a particular 
Hamiltonian according to the approach discussed in \cite{PKtmfVar}, 
\cite{ChKul}.

The properties of the generalized Cremmer-Gervais $R$-matrix
$R(q,p,\nu)$ and the corresponding quantum group ${\cal A}(R)$
and quantum algebra (a dual to the  Hopf algebra ${\cal A}(R)$)
are discussed in Sec.2. The quantum group ${\cal A}(R)$ 
and quantum algebra 
corresponding to the original $R_{CG} = R(q,q^{1/3},\nu) $ 
were studied already in \cite {Balog}. The next section deals with the
$q$-oscillator interpretation of the quantum space
${\bf C}_q^3(R(q,p,\nu)).$ The different values of the
deformation ($q$) and twist ($p,\nu$) parameters correspond to 
the different choice of 
the $q$-oscillator generators \cite {KulD,DKulIJMP}. An integrable
spin-chain Hamiltonian is discussed in Sec.4. Although this
Hamiltonian is not Hermitian due to the inequality $({\cal P}
R_{CG})^{\dagger } \neq {\cal P}R_{CG},$ its structure 
coincides with the one of \cite{Ritten}
where non-Hermitian Hamiltonians were used 
to describe some kinetic processes. Let us mention 
also a couple of recent papers \cite {StoK,MS,Fio} 
demonstrating an active interest in different applications 
of the Drinfeld twist. 

{\bf Acknowledgements.} PPK thanks the Helsinki Institute of 
Physics for the warm hospitality and support. This research was 
supported in part by the RFFI grants N 96-01-00851 (PPK) and N 
97-01-01152 (EVD). 

\section{Properties ~of ~the ~generalized 
~Cremmer-Gervais ~$R$-matrix}
\setcounter{equation}{0}

Although in the framework of the quantum group 
theory one uses such general
objects as universal $R$-matrix ${\cal R}\in {\cal A}\otimes 
{\cal A},$ where ${\cal A}$ is a quasitriangular Hopf algebra 
\cite{DrQG} and universal twist
 element ${\cal F}\in {\cal A}\otimes {\cal A},$ with 
appropriate properties (see, e.g. \cite{DR-Tw,ReshTw}, we will be 
working with finite dimensional 
matrices which are images of ${\cal R}$ and ${\cal 
F}$ in corresponding representations $\rho ,\pi $ of the 
quasi-triangular Hopf algebra ${\cal A}$

\begin{equation}
R=(\rho \otimes \pi ){\cal R}\,,\qquad F=(\rho \otimes \pi ){\cal F}\,. 
\label{a2.1}
\end{equation}

A twisted $R$-matrix ${\cal R}^{({\cal F})}$ is given by the 
twist transformation \cite{DR-Tw,ReshTw} 
\begin{equation} 
{\cal R}^{({\cal F})}={\cal F}_{21}{\cal RF}^{-1}  
\label{a2.2} 
\end{equation}
of the original $R$-matrix ${\cal R},$ where ${\cal 
F}_{21}={\cal PFP},$ and ${\cal P}$ is the permutation map in 
${\cal A}\otimes {\cal A}$. The twist element satisfies the 
relations  in ${\cal A}\otimes {\cal A}$ \cite{DR-Tw,ReshTw} 
$$ 
(\epsilon \otimes id)\, {\cal F} = (id \otimes \epsilon)\, {\cal F} = 1\,, 
$$ 
and in ${\cal A}\otimes {\cal A}\otimes {\cal A}$ 
$$
{\cal F}_{12}\,(\Delta \otimes id)\, {\cal F} = 
{\cal F}_{23}\, (id \otimes \Delta) \,{\cal F} \,. 
$$

In the fundamental representation of the $sl(3)$ algebra the generalized 
Cremmer-Gervais $R$-matrix $R_{CG}$ is defined by the standard 
$R(q)$ (\ref{a1.1}) and the twist matrix $F$ \cite{Jacob}
\begin{equation}
F=\left(
\begin{array}{ccccccccc}
1 & 0 & 0 & 0 & 0 & 0 & 0 & 0 & 0 \\
0 & 1 & 0 & 0 & 0 & 0 & 0 & 0 & 0 \\
0 & 0 & q/p & 0 & p\nu & 0 & 0 & 0 & 0 \\
0 & 0 & 0 & p & 0 & 0 & 0 & 0 & 0 \\
0 & 0 & 0 & 0 & p & 0 & 0 & 0 & 0 \\
0 & 0 & 0 & 0 & 0 & 1 & 0 & 0 & 0 \\
0 & 0 & 0 & 0 & 0 & 0 & p & 0 & 0 \\
0 & 0 & 0 & 0 & 0 & 0 & 0 & p & 0 \\
0 & 0 & 0 & 0 & 0 & 0 & 0 & 0 & 1
\end{array}
\right), \qquad F_{21}=\left(
\begin{array}{ccccccccc}
1 & 0 & 0 & 0 & 0 & 0 & 0 & 0 & 0 \\
0 & p & 0 & 0 & 0 & 0 & 0 & 0 & 0 \\
0 & 0 & p & 0 & 0 & 0 & 0 & 0 & 0 \\
0 & 0 & 0 & 1 & 0 & 0 & 0 & 0 & 0 \\
0 & 0 & 0 & 0 & p & 0 & 0 & 0 & 0 \\
0 & 0 & 0 & 0 & 0 & p & 0 & 0 & 0 \\
0 & 0 & 0 & 0 & p\nu & 0 & q/p & 0 & 0 \\
0 & 0 & 0 & 0 & 0 & 0 & 0 & 1 & 0 \\
0 & 0 & 0 & 0 & 0 & 0 & 0 & 0 & 1
\end{array}
\right),  \label{a2.3}
\end{equation}
depending on two arbitrary parameters $p$ and $\nu$
\begin{equation}
R(q,p,\nu)=F_{21}R(q)F^{-1}.  
\label{a2.4}
\end{equation}

A twist 
transformation of $R$-matrix related to the braid group
${\check R}={\cal P}R$,
\begin{equation}
{\check R}^{(F)}={\cal P}R^{(F)}={\cal P}F_{21}RF^{-1}=
F{\check R}F^{-1} \,, 
\label{a2.5}
\end{equation}
is a similarity transformation. Hence, the spectral
characteristics of ${\check R}$ and ${\check R}^{(F)}$ are the
same. In particular ${\check R}_{CG}$ also satisfies the
Hecke condition $(\omega = q-q^{-1})$
\begin{equation}
\begin{array}{l}
{\check R}^2=I+ \omega {\check R}, \\
{\check R}=qP^{(+)}- q^{-1} P^{(-)},
\end{array}
\label{a2.6}
\end{equation}
where $P^{(\pm )}$ are projectors of $rankP^{(+)}=n(n+1)/2$ and
$rankP^{(-)}=n(n-1)/2,$ $n=3$ for the case of $sl(3),$ we are
interested in.

Using the FRT approach \cite{FRT}, the standard notations of the QISM 
and a particular $R$-matrix, 
one can define a quantum group ${\cal A}(R)$ with a 
matrix $T=\left\{ T_{ij}\right\} $ of the generators satisfying 
the FRT-relation 
\begin{equation}
R_{12}T_1T_2=T_2T_1R_{12,}\quad {\rm or\quad }
\left[ {\check R},T\otimes T\right] {\rm =0}.
\label{a2.7}
\end{equation}
The quantum group ${\cal A}(R)$ associated with the
original Cremmer-Gervais $R$-matrix was defined and discussed
already in \cite{Cr-Ger} and \cite{Balog}. In particular, the
constructed quantum determinant was central. In the case of
$R(q,p,\nu)$ (\ref{a2.4}) the quantum determinant
$\det\nolimits_qT$ is not central. It can be defined by fusion
procedure \cite{KulSk}, projecting $T^{\otimes 3} = 
T_1 T_2 T_3 $ to
one-dimensional space of the $q$-antisymmetrizer
\begin{equation}
P_{123}^{(-)}\simeq P_{12}^{(-)}\left( \left( q+\frac 1q\right)
P_{23}^{(-)}-I\right) P_{12}^{(-)}.  
\label{a2.8}
\end{equation}
The Yang - Baxter equation for $R_{12}$ and the defining 
relations (\ref {a2.7}) mean that the blocks $R_{ik,jl}$ with 
fixed $i,j$ give rise to a representation of the ${\cal A}(R)$ 
generators $T_{ij}.$ The corresponding quantum determinant of 
$R(q,p,\nu)$ is not proportional to the identity $3{\times}3$ matrix 
although it is diagonal
\begin{equation}
\det\nolimits_q R(q,p,\nu)=q \;{\rm diag}
(\frac q{p^3},1,\frac{p^3}q).
\label{a2.9}
\end{equation}
Projecting the commutation relation of $T_4=I^{\otimes 3}\otimes 
T$ with $T^{\otimes 3}=T_1T_2T_3=T\otimes T\otimes T\otimes I$ 
\begin{equation}
\prod\limits_{j=1}^3R_{j4}(T^{\otimes 3}\,T_4) = (T_4\,
T^{\otimes 3})\prod\limits_{j=1}^3R_{j4}\,,   
\label{a2.10}
\end{equation}
 and using the $q$-antisymmetrizer $P_{123}^{(-)}$ (\ref{a2.8}), 
one gets a compact matrix 
form of the commutation relations between $\det\nolimits_qT$ and 
$T_{ij},$ provided the former is invertible

\begin{equation}
\left( \det\nolimits_qT\right) T
\left( \det\nolimits_qT\right) ^{-1}=
\left( \det\nolimits_qR\right) ^{-1}T\left(
\det\nolimits_qR\right) ,
\label{a2.11}
\end{equation}
where the r.h.s. is the usual matrix product of $3{\times}3$ matrices $T$ and
(\ref{a2.9}). 

For discussing a quantum group ${\cal A}(R)$ coaction on quantum
spaces, it is important to have a $*$-operation on
${\cal A}(R)$ and on the  quantum space. 
Taking the deformation parameters real,
$q,p,\nu \in {\bf R}$ and applying a $*$-operation to the
defining FRT-relation ($\ref{a2.7}$) one gets
\begin{equation}
R_{12}T_2^{*}T_1^{*}=T_1^{*}T_2^{*}R_{12},
\label{a2.12}
\end{equation}
where $T^* = \{(T_{ij})^*\}$. By the similarity 
transformation of this relation with matrices
${\cal P}$ and $C\otimes C,$
where $C_{ij}=\delta _{i,4-j},$ the relation ($\ref{a2.12}$)
is reduced to
the original form ($\ref{a2.7}$)
\begin{equation}
R_{12}\widetilde{T}_1^{*}\widetilde{T}_2^{*}=
\widetilde{T}_2^{*}\widetilde{T}_1^{*}R_{12}\,, \quad 
\widetilde{T}^{*}=CT^{*}C \,.
\label{a2.13}
\end{equation}
Hence, the $*$-operation on ${\cal A}(R_{CG})$ can be
defined as follows (cf. \cite {Irac})
\begin{equation}
T^{*}=CTC \,.  
\label{a2.14}
\end{equation} 
This form of the $*$-operation is valid also for 
the quantum group in the $SL(n)$ 
case with the matrix $C_{ij}=\delta _{i,n + 1- j}$. 
The corresponding quantum determinant of the $R$-matrix 
is also diagonal $n{\times}n$ matrix 
\begin{equation}
\det\nolimits_q R(q,p,\nu) \simeq {\rm diag}
(1,\, q^2(p/q)^n,\,.\,.\,.,\, (q^2(p/q)^n)^{n-1}).
\label{a2.15} 
\end{equation} 

Particular relations among parameters $q$ and $p$ 
result in the possibility of further
restrictions on the ${\cal A}(R_{CG})$ generators 
(see the next section). 

The dual to ${\cal A}(R_{CG})$ quasitriangular Hopf algebra 
with  $ p = q^{1/3}$ 
was described in \cite{Balog}. The knowledge of the twist 
element $\cal F$ permits to write the entries of the 
$L$-matrices in terms of the standard $sl_q(3)$ algebra 
generators. These quantities will not be used below and 
we end the discussion of the general properties of the 
quantum group and quantum algebra related to 
the generalized Cremmer-Gervais $R$-matrix $R(q,p,\nu)$.

\section{Quantum ~linear ~space ~of ~${\cal A}(R_{CG})$ ~as 
 ~covariant ~~~~~ deformed oscillator} 
\setcounter{equation}{0}

The quantum linear space of ${\cal A}(R_{CG})$ as covariant
algebra ( ${\cal A}(R_{CG})$-comodule algebra ), according 
to the general 
approach \cite{FRT,WZ}, is generated by three elements
$x_i,\,i=1,2,3;\, x_1 = A,\, x_2 = K,\, x_3 = A^{\dag}$ 
which satisfy relations ($X^t=(A,\, K,\, A^{\dag})$)
\begin{equation}
\label{a3.1}
R(q,p,\nu)X_1X_2=qX_2X_1 = q\,{\cal P}X_1X_2\,,
\end{equation}
where $X_1X_2=X\otimes X$ and $X_2X_1={\cal P}X_1X_2$.

In terms of the generators the defining relations are
\begin{equation}
\label{a3.2}
KA=pqAK,\qquad KA^{\dag}=\frac 1{pq}A^{\dag}K,\qquad
AA^{\dag}-p^{-2}A^{\dag}A=\nu K^2.
\end{equation}

The $*$-operation is: $K^*=K,\quad (A^{\dag})^*=A$ . It is 
consistent with the $*$-operation of the quantum group 
${\cal A}(R_{CG})$-coaction: 
\begin{equation}
\label{a3.2x}
X^{*} = C\,X \,, \quad  T^{*} = C\,T\,C\, \quad  (T\,X)^{*} = C\,(T\,X) \,. 
\end{equation}

Although there are three parameters $q, p, \nu $ in the defining 
relations (3.2), this ${\cal A}(R_{CG})$-comodule algebra is nothing but 
a $q$-oscillator algebra ${\cal A}_q$. Using for the $q$-oscillator algebra 
the Arik-Coon generators (for a convenience of further identification 
we use $q^2$ as the deformation parameter of ${\cal A}_q$) 
\begin{equation}\label{a3.3}
[N, a] = - a\,,\qquad [N, a^{\dag}] = a^{\dag}\,,\qquad
a a^{\dag}-q^{2}a^{\dag} a = 1 
\end{equation}
and doing a one-parameter transformation \cite{KulD} 
\begin{equation}\label{a3.4}
a(\lambda) = q^{-\lambda N} a,\qquad a^{\dag}(\lambda) = a^{\dag} q^{-\lambda N}  
\end{equation}
one gets 
\begin{equation}\label{a3.5}
a(\lambda) a^{\dag}(\lambda) - q^{2( 1 - \lambda )} a^{\dag}(\lambda) a(\lambda) 
=  q^{-2\lambda N}\,. 
\end{equation}

Hence, the relations among generators and parameters of ($\ref{a3.2}$) 
and ($\ref{a3.4}$) are as follows: 
\begin{equation}\label{a3.5x}
A = a(\lambda)\,, \quad A^{\dag} = a^{\dag}(\lambda)\,, \quad  
K^2 = \nu^{-1} q^{-2\lambda N}\,, \quad  p = q^{\lambda - 1}\,. 
\end{equation}
The parameter of the harmonic oscillator deformation is $q$, while 
$p = q^{\lambda - 1}$ refers to a choice of the $q$-oscillator algebra 
generators. Let us point out some special values of these parameters 
which are of particular interest: \\ 
1) $p = 1/q$, $K$ is a central element, the Arik-Coon generators; \\
2) $p = q^{-1/2}, \lambda = 1/2 $ , the Macfarlane-Biedenharn generators ; \\ 
3) $ p = q^{1/3}, \lambda = 4/3 $ , one could call this case 
the Cremmer-Gervais generators, for the quantum determinant of the 
preceding Sec. is central for ${\cal A}(R_{CG})$; \\ 
 4) $q = 1 $ (but $p$ is arbitrary, due to appropriate limit of $\lambda$ ) 
 leads to the standard harmonic oscillator with nonstandard generators. 
 
By construction, these generators ($\ref{a3.2}$) 
of ${\cal A}_q$ are covariant with respect 
to the linear transformation (coaction) of the quantum group 
${\cal A}(R_{CG})$ : $X' = T X$ 
\begin{eqnarray} 
A' &=& T_{11} A + T_{12} K + T_{13}A^{\dag}\,, \nonumber \\
K' &=& T_{21} A + T_{22} K + T_{23}A^{\dag}\,, \\
(A^{\dag})' &=& T_{31} A + T_{32} K + T_{33}A^{\dag}\,. \nonumber 
\label{a3.6}
\end{eqnarray}
 
 This is a $q$-analogue of the Bogoliubov transformation of the canonical 
 commutation relations. However, due to the noncommutativity of the 
 corresponding coefficients $ T_{ij}$ 
 the physical meaning of this transformation 
 seems different from the standard case, where it was used to transform 
 a general quadratic Hamiltonian in question to a simpler, canonical form. 
 The transformation ($\ref{a3.6}$) is different also from the transformation 
 of the $q$-oscillator generators with coefficients depending on the 
 operator $N$ \cite{FI} for in our case $ [T_{ij}\,, X_k] = 0$. 
 
 To illustrate the difference in the interpretation, 
 let us start with the simplest case 4) of the 
 standard harmonic oscillator with nonstandard generators. This case, having 
 trivial deformation parameter $q = 1$, is reduced to a pure twist 
 ${\cal R} = {\cal F}_{21}{\cal F}^{-1} $. However, 
 the generators of the quantum group are non-commuting and one can construct 
 easily irreducible representations of ${\cal A}(R_{CG})$. To have deformation 
 (twist) of the group of the standard Bogoliubov transformations, 
 one can put now two out of nine generators 
 of ${\cal A}(R_{CG})$ to zero $ (T)_{2j} = 0\,, j = 1\,, 3$. Then the quadratic 
 relations of the remaining seven generators can be 
 realized by operators of a Weyl pair  
 $U V = p V U $ with a unitary operator $U$ and a Hermitian $V$. Any irreducible 
 representation of the reduced ${\cal A}(R_{CG})$ is given by these operators and 
 four parameters of the standard Bogoliubov transformation. 
 
 Any Hamiltonian $H(A^{\dag}\,, A)$ (e.g. $H_0 = A^{\dag}\,A$ ) will be mapped 
 under ${\cal A}(R_{CG})$-coaction 
 into an operator in the tensor product space of the harmonic oscillator space 
 and the space of the Weyl pair operators $U\,, V$. The irreducible 
 representations of the complete quantum group ${\cal A}(R_{CG})$ 
 is difficult to construct, 
 because it is noncompact $*$-algebra (for the compact quantum $SU(n)$ 
 group see \cite {Soi}). 
 
 We are using this example to illustrate once more the difficulties with 
 physical interpretation of the group-coaction deformations 
 \cite {PKtmfVar, ChKul}, when one tries to 
 connect this coaction with symmetry properties of the 
 covariant system. This point is especially important while 
 deforming the kinematical groups \cite {AKR}). 

%


\section{Integrable models of the QISM}
\setcounter{equation}{0}

According to the QISM \cite{LDF},\cite{KulSk} one can construct 
a variety of integrable models, related to a given solution of 
the YBE (an $R$-matrix). The importance of the corresponding 
spin-chain models is connected
with possibilities to get different field-theoretical integrable 
models using different limiting procedures \cite{LDF}. The extra 
parameters adding to the $R$-matrix by a twist $F$ are quite
useful in this way.  It was pointed out \cite{StoK}, that the
twist of an $R$-matrix
\[
R(u)\rightarrow
F_{21}R_{12}(u)F^{-1}\equiv R^{(F)}(u)\, , 
\]
preserves the
regularity property \cite{KulSk}: the existence of the spectral
parameter value $u=u_0,$ where $R(u_0)$ is proportional to the
permutation operator ${\cal P}.$ It is obvious, that from
$R(u_0)$ $\simeq {\cal P}$ it follows
\[
R^{(F)}(u_0)=F_{21}R_{12}(u_0)F^{-1}\simeq {\cal P}.
\]
This regularity property is important to get local integrals of 
motion (see e.g. [1,2,3]).

 The Hecke condition leads to the spectral parameter dependent 
solution of the YBE through the Yang-Baxterization procedure 
\cite{Jones}  
\begin{equation} 
{\check R}(u)=u{\check R}-\frac
1u{\check R}^{-1}=(u-u^{-1}){\check R}_{CG}+ \omega u^{-1}I.
\end{equation}
Let us remind that
${\check R}(u)$ satisfies the YBE in the braid group
form.

The $L$-operator of the integrable spin chain coincides with the
$R$-matrix $R(u)$ and the density of the Hamiltonian is
($u_0=\pm 1$)
\[
H=\sum_nh_{n,n+1},
\]
\[ \begin{array}{c}
h_{1,2}\simeq {\cal P}\frac d{du}R(u)|_{u=u_0}=
2 {\check R}_{CG}+const. \\[8pt]
 \simeq \left(
\begin{array}{ccccccccc}
q & 0 & 0 & 0 & 0 & 0 & 0 & 0 & 0 \\
0 & 0 & 0 & 1/p & 0 & 0 & 0 & 0 & 0 \\
0 & 0 & 0 & 0 & q\nu & 0 & q/p^2 & 0 & 0 \\
0 & p & 0 & \omega & 0 & 0 & 0 & 0 & 0 \\
0 & 0 & 0 & 0 & q & 0 & 0 & 0 & 0 \\
0 & 0 & 0 & 0 & 0 & 0 & 0 & 1/p & 0 \\
0 & 0 & p^2 /q & 0 & -p^2\nu/q & 0 &
\omega & 0 & 0 \\
0 & 0 & 0 & 0 & 0 & p & 0 & \omega & 0 \\
0 & 0 & 0 & 0 & 0 & 0 & 0 & 0 & q
\end{array}
\right) \\[8pt]
\begin{array}{cl} 
=&\sum\limits_{i=1}^3 q\,e_{ii}
\otimes e_{ii}+\left( p\,e_{21}\otimes e_{12} 
+p^{-1} e_{12}\otimes e_{21}\right) +
q \nu \,\left( p\,e_{12}\otimes e_{32}-\left( p/q\right)
^2e_{32}\otimes e_{12}\right) \\[5pt]
&+ \left( p^2/q\,e_{31} \otimes e_{13} 
+q/p^{2} e_{13}\otimes e_{31}\right)
+ \omega \sum\limits_{a<b}e_{bb}\otimes e_{aa} \,.
\end{array}
\end{array}
\]

Due to the fact that the twist transformation leads to the explicit 
expression for the new $R$-matrix in terms of the old one and the twist 
element $\cal F$, one can get connection between quantum scattering data 
(the transition matrix) of initial and deformed integrable models 
\begin{equation}\label{a4.4}
 T_{t}(u) = (\rho \otimes \pi) (id \otimes \Delta^{(N)}_t) 
 {\cal F}_{21} {\cal R} {\cal F}^{-1}\,, 
\end{equation}
where $ \rho $ and $\pi$ are representations of the quantum algebra in the 
auxiliary and quantum space respectively  \cite{StoK}. However, the general 
explicit expression for $ T_{t}(u)$ looks rather cumbersome. The formulae (4.2) 
includes $(N-1)$ iterations of the twisted coproduct $ \Delta_t = {\cal F} 
\Delta \,{\cal F}^{-1} $ . Introducing notations 
\begin{equation}\label{4.4}
(id \otimes \Delta) {\cal F} = {\cal F}_{1,23}\,, \quad
(id \otimes \Delta^{(3)}) {\cal F} = (id \otimes id \otimes \Delta) 
(id \otimes \Delta){\cal F} = {\cal F}_{1,234} \quad etc. \,, 
\end{equation}
where the index 1 refers to the auxiliary quantum algebra (auxiliary space 
with irrep $\rho $). Using these notations one gets for the iteration of the  
twisted coproduct $ \Delta_t $ 
\begin{equation}\label{4.5}
 \Delta_t^{(3)} =  (id \otimes \Delta_t) \Delta_t = 
 {\cal F}_{23} (id \otimes \Delta)({\cal F} \Delta {\cal F}^{-1}) 
 {\cal F}_{23}^{-1} = {\cal F}_{23} {\cal F}_{1,23} \Delta^{(3)} 
{\cal F}_{1,23}^{-1} {\cal F}_{23}^{-1} = X_{123} \Delta^{(3)} (X_{123})^{-1}\,. 
\end{equation} 
Hence, the formulae for the transition matrix $T_{t}(u)$ of the twisted 
(deformed) spin chain on $(N - 1)$ sites is the following 
\begin{equation}\label{4.6}
 T_{t}(u) = (\rho \otimes \pi^{\otimes (N - 1)})\, X_{12...N} 
((id \otimes \Delta^{(N)}) {\cal F}_{21} {\cal R} {\cal F}^{-1}) X_{12...N}^{-1}\,.
\end{equation} 

Although the iterated coproduct is a similarity transformation of the 
original one (4.4), the transition matrices are related by a more 
complicated transformation, since the factors in (4.5) are not the inverse 
of each other. The Yang-Baxter algebra of the quantum scattering data 
$( T_{t}(u) )_{ij}$ is more complicated than the $SL(3)$-spin chain \cite {KR} 
due to the extra non-zero elements of $R_{CG}(u)$ proportional to $\nu$. 
However, as in the case of the deformed (twisted) $XXX_{\xi}$-spin chain 
\cite {Ritten, StoK}, the spectrum of the transfer matrix $t(u) = tr T_{t}(u)$  
is expected to be the same as for the $SL(3)$-spin chain, 
taking into account the 
changed eigenvalues of the diagonal elements of $ T_{t}(u)$ on the 
reference (vacuum) state $(0,\, 0,\, 1)^t$. 
The Bethe equations defining the two sets of quasimomenta  have  
similar structure with obvious changes due to the parameter $p$ 
\cite{KR}. The second 
parameter $\nu$ enters into the eigenvectors and adjoint vectors. 

\section{Conclusion} 
\setcounter{equation}{0}

The twist transformation in the theory of quantum groups preserves 
the algebraic sector of the twisted Hopf algebra. Hence the representations 
of the twisted algebra are the same. However, the 
change of the coproduct has deep consequences for the R-matrix 
structure and  for the corresponding integrable models, which rely on 
the tensor product of these representations. The generalized 
Cremmer-Gervais $R$-matrix, being a twist of the standard one, 
leads to more complicated integrable systems. Two of these systems 
were studied in this paper: i) the covariant $q$-oscillator algebra,  
and ii) an integrable $sl(3)$ spin chain. The covariance of the 
$q$-oscillator algebra under the $q$-Bogoliubov transformation 
results in preservation of the ${\cal A}_q$ structure in an {\it extended}  
dynamical system. We hope that further research and explicit 
knowledge of the twist element permits to calculate most of the 
characteristics of the twisted models in terms of the initial ones.

\end{document}